\documentclass[12pt]{article}
\usepackage{amssymb,amsthm,amsmath,amsfonts,latexsym,tikz,hyperref,ytableau,rotating,authblk}
\usepackage[hmargin=.8in,vmargin=1in]{geometry}
\usepackage{diagbox}

\newtheorem{thm}{Theorem}[section]
\newtheorem{prop}[thm]{Proposition}
\newtheorem{cor}[thm]{Corollary}
\newtheorem{lem}[thm]{Lemma}
\newtheorem{conj}[thm]{Conjecture}
\newtheorem{exa}[thm]{Example}

\DeclareMathOperator{\lwi}{lwi}
\DeclareMathOperator{\Par}{Par}
\DeclareMathOperator{\RC}{RC}
\DeclareMathOperator{\rw}{rw}

\newcommand{\fP}{\mathfrak P}

\newcommand{\ben}{\begin{enumerate}}
\newcommand{\een}{\end{enumerate}}
\newcommand{\ble}{\begin{lem}}
\newcommand{\ele}{\end{lem}}
\newcommand{\bth}{\begin{thm}}
\renewcommand{\eth}{\end{thm}}
\newcommand{\bpr}{\begin{prop}}
\newcommand{\epr}{\end{prop}}
\newcommand{\bco}{\begin{cor}}
\newcommand{\eco}{\end{cor}}
\newcommand{\bcon}{\begin{conj}}
\newcommand{\econ}{\end{conj}}
\newcommand{\bde}{\begin{defn}}
\newcommand{\ede}{\end{defn}}
\newcommand{\bex}{\begin{exa}}
\newcommand{\eex}{\end{exa}}
\newcommand{\barr}{\begin{array}}
\newcommand{\earr}{\end{array}}
\newcommand{\btab}{\begin{tabular}}
\newcommand{\etab}{\end{tabular}}
\newcommand{\beq}{\begin{equation}}
\newcommand{\eeq}{\end{equation}}
\newcommand{\bea}{\begin{eqnarray*}}
\newcommand{\eea}{\end{eqnarray*}}
\newcommand{\bal}{\begin{align*}}
\newcommand{\bce}{\begin{center}}
\newcommand{\ece}{\end{center}}
\newcommand{\bpi}{\begin{picture}}
\newcommand{\epi}{\end{picture}}
\newcommand{\bpp}{\begin{picture}}
\newcommand{\epp}{\end{picture}}
\newcommand{\bfi}{\begin{figure} \begin{center}}
\newcommand{\efi}{\end{center} \end{figure}}
\newcommand{\bprf}{\begin{proof}}
\newcommand{\eprf}{\end{proof}\medskip}

\newcommand{\bsl}{\begin{slide}{}}
\newcommand{\esl}{\end{slide}}
\newcommand{\bfr}{\begin{frame}}
\newcommand{\efr}{\end{frame}}

\newcommand{\eqqed}[1]{$\rule{1ex}{0ex}\hfill{\dil#1}\hfill\qed$}

\newcommand{\hso}[1]{\hspace{-1pt}}
\newcommand{\vs}[1]{\vspace{#1}}

\newcommand{\sbe}{\subseteq}

\newcommand{\ptn}{\vdash}

\newcommand{\lte}{\unlhd}

\newcommand{\case}[4]{\left\{\barr{ll}#1&\mbox{#2}\\#3&\mbox{#4}\earr\right.}
\newcommand{\fl}[1]{\lfloor #1 \rfloor}
\newcommand{\ce}[1]{\lceil #1 \rceil}

\def\<{\langle}
\def\>{\rangle}

\newcommand{\ra}{\rightarrow}

\newcommand{\al}{\alpha}
\newcommand{\be}{\beta}

\newcommand{\de}{\delta}
\newcommand{\ep}{\epsilon}

\newcommand{\la}{\lambda}

\newcommand{\1}{{\bf 1}}

\newcommand{\bbN}{{\mathbb N}}
\newcommand{\bbP}{{\mathbb P}}

\newcommand{\cL}{{\cal L}}

\newcommand{\cP}{{\cal P}}

\newcommand{\fS}{{\mathfrak S}}

\DeclareMathOperator{\des}{des}

\newcommand{\dil}{\displaystyle}

\begin{document}
\pagestyle{plain}

\title{Centralizers in the plactic monoid
}
\author[1]{Bruce E. Sagan}
\author[2]{Alexander N. Wilson}
\affil[1]{Department of Mathematics, Michigan State University, East Lansing, MI 48824}
\affil[2]{Department of Mathematics, Oberlin College, Oberlin, OH 44074}

\date{\today\\[10pt]
	\begin{flushleft}
	\small Key Words: centralizer, jeu-de-taquin, Knuth equivalence, partition, plactic monoid, Robinson-Schensted-Knuth correspondence, 
	                                       \\[5pt]
	\small AMS subject classification (2020):  05E99  (Primary) 05A05, 05A15 (Secondary)
%https://mathscinet.ams.org/mathscinet/msc/msc2020.html
	\end{flushleft}}

\maketitle

\begin{abstract}
Let $u$ be a word over the positive integers.  Motivated in part by a question from representation theory, we study the centralizer set of $u$ which is
$$
C(u) = \{w \mid \text{$uw$ is Knuth-equivalent to $wu$}\}.
$$
In particular, we give various necessary conditions for $w$ to be in $C(u)$.  We also characterize $C(u)$ when $u$ has few letters, when it has a single repeated entry, or when it is a certain type of decreasing sequence.
We consider $c_{n,m}(u)$, the number of $w\in C(u)$ of length $n$ with $\max w\le m$.
We prove that for $|u|=1$ the value of this function depends only on the relative sizes of $u$ and $m$ and not on their actual values.  And for various $u$ we use Stanley's theory of poset partitions to show that, for fixed $n$,  $c_{n,m}(u)$ is a polynomial in $m$ with certain degree and leading coefficient.  We end with various conjectures and directions for further research.
\end{abstract}

%%%%%%%%%%%%%%%%%%%%%%%%%%%%%%%%
\section{Introduction}	
%%%%%%%%%%%%%%%%%%%%%%%%%%%%%%%%

Let $\bbP =\{1,2,3,\ldots\}$ and $\bbN=\bbP\uplus\{0\}$ denote the positive and nonnegative integers, respectively.  For $n\in\bbN$ we let
$$
[n] = \{1,2,\ldots,n\}.
$$
In addition, for any set $S$ we will use either $\#S$ or $|S|$ to denote the cardinality of $S$.
We apply the same notation to words $w$ over $S$ and call $|w|$ the {\em length} of $w$.
Finally, we let $S^*$ be the {\em Kleene closure} of $S$, that is, all words with elements from $S$.

We will assume the reader is familiar with the Robinson-Schensted-Knuth (RSK) correspondence as well as Sch\"utzenberger's jeu-de-taquin (jdt).  Background on these operations can be found in the texts of Sagan~\cite{sag:sg,sag:aoc} or Stanley~\cite{sta:ec2}.  In particular, if $w\in\bbP^*$ then we will let 
$P(w)$ denote the insertion tableau of $w$ under RSK.  Recall that $v,w\in\bbP^*$ differ by a {\em Knuth transposition} if there are words $x,y$ and elements $a,b,c$ such that either
$$
\text{$v = xacby$ and $w=xcaby$ with $a\le b<c$},
$$
or
$$
\text{$v = xbacy$ and $w=xbcay$ with $a< b\le c$}.
$$
Furthermore, we say that $v,w$ are {\em Knuth-equivalent}, written $v\equiv w$, if one can obtain $w$ from $v$ by applying a sequence of Knuth transpositions.  When Knuth introduced this equivalence relation~\cite{knu:pmg}, he proved that
$$
v\equiv w \text{ if and only if } P(v) = P(w).
$$
The {\em plactic monoid} is $\bbP^*$ modulo Knuth equivalence.
It was first considered from this perspective by Lascoux and Sch\"utzenberger~\cite{LS:mp}.

Given a word $u\in\bbP^*$, our primary object of study will be the {\em centralizer of $u$} in the plactic monoid  which is
$$
C(u)=\{w \mid \text{$uw\equiv wu$}\},
$$
or equivalently
$$
C(u)=\{w \mid P(uw)=P(wu)\}.
$$
In particular, we wish to characterize  $C(u)$ for certain $u$ and also consider the enumerative properties of the integers 
\beq
\label{cnmu}
c_{n,m}(u) =\#\{w\in C(u) \mid \# w = n \text{ and } \max w\le m\}.
\eeq
Beside the fact that $C(u)$ is a natural set to study, our research  
is motivated by work in preparation by the second author and Nate Harman concerning commuting crystal structures on ``lexicographic bitableaux,'' semistandard tableaux filled with entries in $[m]\times[n]$ ordered lexicographically. In this setting, it is natural for one crystal operation to transform a reading word of the tableau by cutting out a subword and pasting it in a different location. In order for these crystals to commute, the transformed word must be Knuth-equivalent to the original reading word.

The rest of this paper is organized as follows.  In the next section we will collect some necessary conditions for $w$ to be in $C(u)$ which will prove useful in the sequel.  In section~\ref{csu} we will characterize the $w\in C(u)$ for certain $u$ with  $\#u\le 3$.  In particular, we will describe $C(u)$ for any $u$ of length $1$. 
Next, we will describe $C(u)$ for certain special $u$ of arbitrary length such as those which consist of a single repeated integer or are of the form $m(m-1)\ldots 1$ for some $m\in\bbP$.  Section~\ref{e} is devoted to the study of the $c_{n,m}(u)$ as defined in~\eqref{cnmu}.  In particular, if $|u|=1$ we show that their values depend only on the relative sizes of $m$ and $u$.  Furthermore, we use Stanley's theory of poset partitions to prove that for certain $u$ and fixed $n$, they are polynomials in $m$.  We end with a section containing open problems and conjectures.

%%%%%%%%%%%%%%%%%%%%%%%%%%%%%%%%
\section{Necessary Conditions}
%%%%%%%%%%%%%%%%%%%%%%%%%%%%%%%%

In this section, we collect results giving general constraints on the tableaux $P=P(w)$ for $w\in C(u)$.   In particular, we will give a criterion which will permit us to bound the size of the elements in the first few rows of $P$ by the maximum value in $u$.  Our principle tool here and going forward will be to compare the computation of $P(wu)$ using RSK with the computation of $P(uw)$ using jdt.
In the former, the elements of $u$ are inserted into $P(w)$ using the usual RSK bumping procedure.  In the latter, a skew tableau is formed with $P(u)$ in the southwest and $P(w)$ in the northeast.  The tableau is then brought to left-justified shape using jdt slides.

\begin{figure}
$$
P =\  \begin{ytableau}
    1 & 1 & 1 & 3 & 4 & 4\\
    2 & 3 & 3 & 4 & 5\\
    3 & 5 & 5\\
    4
\end{ytableau}
$$
    \caption{A semistandard Young tableau (SSYT), $P$}
    \label{P}
\end{figure}

Given any sequence $R$ and any element $a$ we let
$$
m_a(R) = \text{ the multiplicity of $a$ in $R$}.
$$
Also, for a semistandard Young tableau (SSYT) $P$ with rows $R_1,R_2,\ldots$, we consider the weak composition 
$$
\al_a(P) = (m_a(R_1),\ m_a(R_2),\ \ldots).
$$
For example, if $P$ is the tableau in Figure~\ref{P},
then
$$
\al_4(P) = (2,1,0,1,0,0,\ldots).
$$
We will compare weak compositions $\al=(\al_1,\al_2,\ldots)$  and $\be=(\be_1,\be_2,\ldots)$ using {\em dominance order} where $\al\preceq\be$ if
$$
\al_1 + \al_2 +\cdots+\al_i\le\be_1 + \be_2 + \cdots+\be_i
$$
for all $i\ge1$.

\ble
    Let $a\neq b$ be distinct positive integers and let $w\in\bbP^*$. Then 
$$
        \alpha_b(P(wa))\preceq\alpha_b(P(w))\preceq\alpha_b(P(aw)).
$$  
\ele

\bprf
    Since $a\neq b$, the number of $b$'s will not change in passing from $P(w)$ to $P(wa)$.  Also, a $b$ can only be bumped to the next row by the RSK algorithm. 
 So, $\alpha_b(P(wa))$ is obtained from $\alpha_b(P(w))$ by subtracting one from the $i$th entry and adding one to the $(i+1)$st entry if  $b$ is bumped from row $i$ to row $i+1$, hence $\alpha_b(P(wa))\preceq\alpha_b(P(w))$.

    For the second inequality, note that  values $b$ can only be slid left or up by the jdt algorithm. It follows that  $\alpha_b(P(aw))$ is obtained from $\alpha_b(P(w))$ by adding one to the $i$th entry and subtracting one from the $(i+1)$st entry if a  $b$ is slid from row $i+1$ to row $i$, hence $\alpha_b(P(aw))\succeq\alpha_b(P(w))$.
\eprf

We can now prove our first necessary condition for when $w\in C(u)$.

\bco
    \label{noBumping}
    
    If $w\in C(u)$ and $b\not\in u$ then 
    $$
        \alpha_b(P(wu))=\alpha_b(P(w))=\alpha_b(P(uw)).
    $$
    In particular, no $b\not\in u$ can be bumped by the insertion of $u$ into $P(w)$ to form $P(wu)$.   And such an element $b$ cannot slide between two rows in the computation of $P(uw)$ by jdt.
\eco
\bprf
Since $w\in C(u)$ we have $P(wu)=P(uw)$.  So, $\al_b(P(wu))=\al_b(P(uw))$ for any $b$.
Now using the fact that $b\not\in u$ and repeated application of the previous lemma finishes the proof.
\eprf

We can now bound the size of certain elements in $P(w)$ for $w\in C(u)$ in terms of the maximum value in $u$.

\ble
\label{max<=m}
    Given $u$ and $w\in C(u)$ we let $P=P(w)$ have rows $R_i$ for $i\geq1$. Also let $m=\max u$. If $u$ contains a subsequence $m, m-1, \dots, m-k+1$, then
$$
    \max R_i\leq m
$$ for $1\leq i\leq k$.
\ele

\bprf
First we claim that if a semistandard tableau $T$ contains an $a$ in a higher row than an $a+1$, then this will continue to be the case after any insertion into $T$.  This is clear if the $a$ and $a+1$ are at least two rows apart since an element can only be bumped from one row to the next.  The other case is if $a$ is bumped from the row directly above $a+1$.  But then the $a$ must bump one of the  elements equal to $a+1$ in the row it enters and so the claim still holds.

Using the assumed subsequence of $u$ and an argument like that in the previous paragraph, we see that in forming $P(wu)$ from $P$ we must have the elements $m,m-1,\ldots,m-k+1$ from $u$ in separate rows with $m$ in the lowest row.  This means $m$ must have traveled through at least the first $k$ rows to its present position.  But if one of these rows contains an element from $P$ larger then $m$, then $m$ would bump the smallest such element to the next row.
This contradicts Corollary~\ref{noBumping} which completes the demonstration.
\eprf

%%%%%%%%%%%%%%%%%%%%%%%%%%%%%%%%
\section{Commuting with short $u$}
\label{csu}
%%%%%%%%%%%%%%%%%%%%%%%%%%%%%%%%

Given a row $R$ of a tableau and a condition $I$ on integers we let
$$
R(I) = \text{ multiset of elements of $R$ satisfying $I$}.
$$
For example, if $u$ is an integer then $R(\le u)$ would be all elements of $R$ which are at most $u$.  More specifically, if $R$ is the second row of the the SSYT, $P$, in Figure~\ref{P} then
$$
R(\le 3) = \{\{2,3,3\}\}.
$$
We will say that cell $(i,j)$ in row $i$ and $j$ is {\em adjacent} to the cells $(i,j+1)$  and $(i+1,j)$ (that is, those which could be next in a jeu-de-taquin path) and similarly with the elements of a tableau in those cells.
\bth
\label{|u|=1:row}
Suppose $u$ consists of a single integer which we also denote by $u$. 
Also, use  $R_1, R_2,\ldots,R_l$ to denote the rows of $P=P(w)$.
Then the  set $C(u)$ is all $w$ such that $P=P(w)$ satisfies
\ben
\item[(a)] $\max R_1\le u$, and
\item[(b)] for $i\ge 1$ we have
$$
\# R_i(<u) = \# R_{i+1}(\le u).
$$
\een
\eth
\bprf
  We first prove that if $P$ satisfies the given restrictions then 
$w\in C(u)$.  That is, we need to prove  $P(wu)=P(uw)$.  
But by (a), $P(wu)$ is obtained from $P$ by appending $u$ to the first row and so has rows $R_1u,R_2,\ldots,R_l$.

To compute $P(uw)$, we perform jeu-de-taquin on the skew tableau with $u$ in the $(l+1,1)$ cell and $P$ in the first $l$ rows starting in column $2$.  
First consider the slide into the cell $(i,1)$ for $i\ge u$.  By reverse induction on $i$, we can assume row $i+1$ is $uR_{i+1}$ so that there is a $u$ in cell $(i+1,1$).
And, by the bound on $i$, the smallest entry in row $i$ is at least $u$ since it comes from  $R_i$ of the semistandard tableau $P$.  It follows that the $u$ in $(i+1,1)$ will move into $(i,1)$.  We now claim that the rest of row $i+1$ will slide left one cell.  This is because, again by reverse induction, $R_{i+2}$ slid one cell left when filling $(i+1,1)$.  Thus the element $x$ below any hole created in row $i+1$ will always be greater than the element $y$ to the hole's right since $x$ was originally directly below $y$ in a semistandard tableau.  This forces a horizontal slide.

Now consider what happens when filling $(u-1,1)$.  Note that row $u$ is currently $uR_u$.
So the slide will start by moving all the elements of $R_{u-1}$ which are smaller than $u$ one cell left.
But by (b) with $i=u-1$, the number of such elements is one less than the number of elements less than or equal to $u$ in $uR_u$.  It follows that the hole created will have a $u$ below and an element at least as large as $u$ to its left.  Thus that hole will be filled by moving the $u$ below into row $u-1$.  Now one can see by an argument in the previous paragraph, that the rest of row $u$ will slide one unit left.  So at the end of this slide row $u$ will be $R_u$ while row $u-1$ will be $R_{u-1}$ with an extra $u$ added.  This is the base case of a reverse induction showing that after sliding into $(i,1)$ for $i<u$, the resulting $(i+1)$st row is $R_{i+1}$ and the $i$th is $R_i$ with a $u$ added.  We omit the induction step which is much like the base case.  It follows that at the end, the second through $l$th rows of the jeu-de-taquin tableau are $R_2,\ldots,R_l$.  And since by (a) all elements of $R_1$ are at most $u$, the result of adding $u$ to $R_1$ must be $R_1u$.  This coincides with the description of $P(wu)$ in the first paragraph, so we are done with this direction.

For the converse, we assume $w\in C(u)$ so that $P(uw)=P(wu)$.
Applying Lemma~\ref{max<=m} with $m=u$ and $k=1$ immediately gives condition (a).
It follows that  
 $P(wu)$, and thus $P(uw)$, has rows $R_1u,R_2,\ldots,R_k$.  
 In computing $P(uw)$ using jeu-de-taquin on the skew tableau, we start out with rows $1$ through $l$ as in $P$ and an extra $u$ in row $l+1$.  Furthermore, in any slide and for any fixed value $a$, at most one $a$ moves between adjacent rows.
Thus for the extra $u$ in row $l+1$ to become an extra $u$ in row $1$ it must be that a $u$ moves from row $i+1$ to row $i$ in the slide filling $(1,i)$ for all $i=l,l-1,\ldots,1$.  It can be shown that this implies (b) using similar ideas as in the first half of the proof.  We leave the details to the reader.
\eprf

We will now give a result which tests for $w$ being in  $C(u)$ when $|u|=1$ by looking at the columns of $P(w)$ rather than the rows.  We use the convention that if  $P$ is an SSYT and $(i,j)$ is a cell  then
$$
P_{i,j} = \text{ element of $P$ in cell $(i,j)$}.
$$
\bth
\label{|u|=1:col}
If $|u|=1$ then
$$
C(u) = \{w \mid \text{every column of $P=P(w)$ contains a $u$}\}.
$$
\eth
\bprf
It suffices to show that $w$ satisfies the column condition above if and only if it satisfies the row conditions (a) and (b) of the previous theorem.  

Assume that (a) and (b) hold and, towards a contradiction, that column $j$ of $P$ does not contain a $u$.  By (a), we must have $P_{1,j}<u$.  Let $i$ be maximal such that $P_{i,j}<u$.  The rows of $P$ being weakly increasing implies that
$\#R_i(<u)\ge j$.  But by the maximality of $i$ and the fact that column $j$ contains no $u$ we also have $\#R_{i+1}(\le u)< j$.  So, using (b),
$$
j \le \#R_i(<u) = \#R_{i+1}(\le u)< j
$$
which is a contradiction.

For the converse, assume every column contains a $u$.  To verify (a) we again argue by contradiction and assume that  $P_{1,j}>u$ for some $j$.  But then strict increase of the columns implies that column $j$ has no $u$ which cannot be.
For (b), consider $R_i$.  If every element of $R_i$ is at least $u$ then every element of $R_{i+1}$ is at least $u+1$ and so both sides of (b) are zero.  Otherwise, let $j$ be maximal such that $P_{i,j}<u$ which implies $\#R_i(<u)=j$.  Since column $j$ must contain a $u$, we have that 
$P_{i+1,j}\le u$ which forces $\#R_{i+1}(\le u)\ge j$.  But if $\#R_{i+1}(\le u)> j$
then we must have  $P_{i+1,j+1}\le u$.  By column strictness $P_{i,j+1}<u$ which contradicts the maximality of $j$.  Thus
$$
\#R_{i+1}(\le u)= j=\#R_i(<u)
$$
as desired.
\eprf

We will now give a characterization of the $w\in C(1)$ which depends directly on $w$ without having to compute $P(w)$. This will also permit us to make a connection with Yamanouchi words. To state these  results, we will need some notation and definitions.  Let
$$
\lwi(w) = \text{ longest length of a weakly increasing subsequence of $w$},
$$
and for $a\in\bbP$
$$
\lwi(w,a) =  \text{ longest length of a weakly increasing subsequence of $w$ of the form $va$}.
$$
For example, if $w=162724534$ then
$\lwi(w)= 5$ because of, for example, the subsequence $12244$ among others.
Also $\lwi(w,3)=4$ as witnessed by $1223$.
Note that to compute $\lwi(w,a)$ one needs only to know the length of a longest weakly increasing subsequence
ending at the rightmost $a$ in $w$:  a weakly increasing sequence ending at an $a$ further to the left can have its length increased by concatenating with the last $a$.
Finally, recall that a word $w$ is {\em Yamanouchi} if every suffix of $w$ has at least as many $i$'s as $(i+1)$'s for all $i\geq1$.

For the next result, we keep the notation in the statement of Theorem~\ref{|u|=1:row}. 
\bco
\label{C(1)}
The following are equivalent.
\ben
\item[(a)] $w\in C(1)$.
\item[(b)] The entries of $R_1$ are all $1$'s.
\item[(c)] $\lwi(w)=\lwi(w,1)$
\een
Furthermore, the following are equivalent.
\ben
\item[(d)] $w\in C(1)\cap[2]^n$.
\item[(e)] $w\in[2]^n$ is Yamanouchi.
\een
\eco
\bprf
(a) $\iff$ (b).  It suffices to show that, when $u=1$, (b) is equivalent to the two conditions in Theorem~\ref{|u|=1:row}.   But the current (b) is Theorem~\ref{|u|=1:row} (a).  And when $u=1$, Theorem~\ref{|u|=1:row} (b) is vacuous since both sides are zero. 

\medskip

(b) $\iff$ (c).  We will prove the forward direction as the reverse is similar.  
By Schendsted's theorem on increasing and decreasing subsequences, $\lwi(w)=\# R_1$.  In fact, if $w_i$ is inserted in column $j$ of $R_1$ during RSK, then $j$ is the longest length of a weakly increasing subsequence of $w$ ending at $w_i$.
By (b), the rightmost $1$ in $w$ is inserted in column $\#R_1$ of the first row so that 
$\lwi(w,1)=\# R_1$.  Combining this with the previous equality finishes the proof.

\medskip

(d) $\iff$ (e).
Since (a) is equivalent to (b), we have that (d) is the same as saying that $R_1$ is all $1$'s and $R_2$ is all $2$'s.  But this means the every $2$ which is inserted during RSK must be bumped into the second row by some following $1$. Thus there is a matching of each $2$ in $w$ with a $1$ which comes afterwards.
 This is the same thing as saying that $w$ is Yamanouchi.
\eprf

For words $u$ of length two or three, we will concentrate on the case when
$u$ consists of $1$'s and $2$'s.  If $\#u=2$ then  $u=11$, $21$, and $22$ will be taken care of by more general results in the next section.  So we will content ourselves with a column characterization for $C(12)$.  To state it we define a {\em singleton $a$-column} to be a column of length $1$ whose entry is $a$.
\bth
\label{C(12):col}
We have that $w\in C(12)$ if and only  if all columns $C$ of $P(w)$ satisfy the following two conditions.
\ben
\item[(a)] If there is a singleton column, $C$,  then $C$ is a singleton $1$-column or a singleton $2$-column, and both types of columns must exist.
\item[(b)] If $\#C\ge2$ then $C$ must contain both $1$ and $2$.
\een
\eth
\bprf
We first show that (a)  and (b) imply $w\in C(12)$.  There are two cases depending on whether $P(w)$ has columns of length one or not.  We will provide details of the former as the latter is similar.

Consider computing $P(w12)$ by RSK.  Conditions (a) and (b) imply that, in the case under consideration, the first row of $P=P(w)$ has both $1$'s and $2$'s and no larger entries.
And from (b), the second row of $P$ contains only of $2$'s.
So insertion of $1$ into $P(w)$ will bump a $2$ into the second row where it will sit at the end of the row.  Then insertion of $2$ into $P(w1)$ will result in the $2$ sitting at the end of the first row.  So in passing from $P$ to $P(wu)$, the multiplicity of $1$ goes up by one in the first row, the multiplictiy of $2$ goes up by one in the second row, and all other multiplicities stay the same.

Now consider computing $P(12w)$ by jdt where $12$ is placed southwest of $P$ which has rows $R_1,R_2,\ldots,R_l$.
We do the slides row by row, starting from the bottom and going up.  It is easy to see that, after filling all rows but the first,
the resulting tableau has rows $R_1$ preceded by blank cells $(1,1)$ and $(1,2)$, $12R_2$,
and all the other rows identical to $P$.  
Note that at this point all rows have the same content as $P$ except row $2$ which has an extra $1$ and an extra $2$.
When cell $(1,2)$ is filled, the cell below contains a $2$ and, by our assumptions, the cell to its left contains a $1$.  So the $1$ and all $1$'s following it will slide left.  Furthermore, since both columns of single $1$'s and single $2$'s exist in $P$ and $R_2$ has only $2$'s, once the $1$'s have slid the resulting empty cell will have no element below it so that the slide will continue to the end of the first row.
Observe that the content of the rows has not changed from before the slide.
Finally we fill the $(1,1)$ cell with the $1$ below it.  Now the rest of the row $2R_2$ will move left since $R_3$ is below it and these were two adjacent rows in a semistandard tableau.  This final change adds a $1$ to the first row and retains the extra $2$ in the second without changing any other multiplicities.  Thus $P(w12)=P(w12)$ and $w\in C(12).$

For the converse, suppose $w\in C(12)$ so that $P(12w)=P(w12)$.  Keeping the notation established in the first part of the proof, from Lemma~\ref{max<=m} 
with $m=2$ and $k=1$
we have that 
$\max R_1\le 2$.  So, any singleton column must consist of a $1$ or $2$.
To show that the rest of (a) as well as (b) holds, we argue by contradiction.
And we just provide details for the former as the latter is similar.
So suppose that $P$ contains a singleton $1$-column but no singleton $2$-column.
Then, by RSK, $P(w12)$ is obtained by appending $12$ to row $R_1$ of $P$.  So the multiplicty of $1$ and $2$ go up by one for $R_1$ and remain the same elsewhere.
But using jdt as before to compute $P(12w)$, the $(1,2)$ cell must be filled by the $1$ to its right.  This causes the whole of $R_1$, to move over one cell and, since there exists at least one column consisting of a single $1$, no elements will be brought up from the second row.  It is now impossible for the second slide to add two elements to $R_1$, a contradiction.  Similarly, assuming that there is a singleton $2$-column but no singleton $1$-column leads to multiplicity problems.
\eprf

We end this section by looking at just one $u$ of length $3$. Since the results and techniques are similar to what we have already seen, we will just sketch the proof.
\bth
\label{C(212):col}
We have that $w\in C(212)$ if and only  if all columns $C$ of $P(w)$ satisfy the following two conditions.
\ben
\item[(a)] All singleton columns are singleton $2$-columnns.
\item[(b)] If $\#C\ge2$ then $C$ must contain both $1$ and $2$.
\een
\eth
\bprf
First assume (a) and (b).  Then computing $P(w212)$ by RSK adds a $1$ and $2$ to the first row, and a $2$ (via bumping) to the second.  It is easy to verify that the computation of $P(212w)$ by jdt will have the same effect, where the lack of any singleton $1$-column will ensure that when filling the $(1,2)$ cell a $2$ from the second row will enter the first to increase the multiplicity as required in comparison to RSK.

For the reverse implication, Lemma~\ref{max<=m} with $m=k=2$ immediately gives condition (b).  And singleton $1$-columns are ruled out by contradiction.
\eprf

%%%%%%%%%%%%%%%%%%%%%%%%%%%%%%%%
\section{Commuting with longer $u$}
\label{clu}
%%%%%%%%%%%%%%%%%%%%%%%%%%%%%%%%

We will now derive characterizations of  $C(u)$ for certain words $u$ of arbitrary length.
We begin by considering the case where $u$ is just the repetition of a single element.  We first have a lemma.  For any word $u$ and $k\ge1$, we let $u^k$ be the concatenation of $k$ copies of $u$.
\ble
\label{u^k}
For any $u\in\bbP^*$ and any $k\in\bbP$ we have
$$
C(u) \sbe C(u^k).
$$
\ele
\bprf
If $w\in C(u)$ then, by induction on $k$,
$$
w u^k \equiv uw u^{k-1} \equiv u u^{k-1} w = u^k w
$$
as desired.
\eprf

We now show that, interestingly, $C(a^k)$ does not depend on $k$ and so can be characterized by the conditions in either Theorem~\ref{|u|=1:row} or Theorem~\ref{|u|=1:col}.
\bth
\label{C(a^k)}
If $a,k\in \bbP$ then
$$
C(a^k) = C(a).
$$
\eth
\bprf
The fact that $C(a)\sbe C(a^k)$ is a special case of the previous lemma.
For the other containment, it suffices to show that if $w\in C(a^k)$ then $w$ satisfies conditions (a) and (b) of Theorem~\ref{|u|=1:row}.  We get (a) immediately from Lemma~\ref{max<=m}.  It follows that, when computing $P(w a^k)$ by RSK, one merely adds $k$ copies of $a$ to the first row of $P=P(w)$.

We now compute $P(a^k w) = P(w a^k)$ using jdt.  The initial skew tableau has $k$ copies of $a$ in a row at the southwest.  And, from the previous paragraph, these must all end up in the first row at the end.  We fill the empty cells starting with the $k$th column from bottom to top, and then working to the left.
From what we have remarked, when working on column $k$ this must bring the $k$th $a$ in $a^k$ from the bottom row to the top.  But, as in the demonstration of Theorem~\ref{|u|=1:row}, this implies that (b) must hold and completes the proof.
\eprf

There is another class of words which have a particularly nice characterization of their centralizers.
\bth
\label{C(m...1)}
We have  $w\in C(m(m-1)\ldots 1)$ if and only if $P=P(w)$ satisfies 
$$
\max R_i \leq m \text{ for all } 1\le i \le m
$$ 
where $R_i$ is the $i$th row of $P$.
\eth
\bprf
Necessity follows immediately from Lemma~\ref{max<=m} with $k=m$.  For sufficiency, let  $u=m(m-1)\ldots 1$ and compute $P(wu)$ by RSK where $P$ satisfies the given restriction.  So, when inserting $m$ into $P$ it will be at least as large as the other elements in $R_1$ and sit at the end of the first row.
Assume, by induction on $j$, that after $m, m-1, \ldots, m-j+1$ are inserted we have a copy of $m-j+i$ in row $i$ for $i\le j$. Then insertion of $m-j$ will cause a cascade of bumps so that, at the end, a copy of each entry has moved down one row and the $m-j$ is now in row $1$.  From this and the assumption on $\max R_i$, we see that the rows of $P(wu)$
are $R_i\uplus\{i\}$ for $i\le m$, and $R_i$ for $i\ge m$.

We now compute $P(uw)$ using jdt.  Note that $P(u)$ is a single column with entries $1,2,\ldots,m$.  It is easy to see that applying the necessary slides merely moves this column up until its elements are added to the first $m$ rows, leaving the rest of the rows invariant.  Thus $P(uw)=P(wu)$ and so $w\in C(u)$ as desired.
\eprf

%%%%%%%%%%%%%%%%%%%%%%%%%%%%%%%%
\section{Enumeration}
\label{e}
%%%%%%%%%%%%%%%%%%%%%%%%%%%%%%%%

We will now use the characterizations derived previously to study the integers $c_{n.m}(u)$ defined by~\eqref{cnmu} which count the number of $w\in C(u)$ of length $n$ with maximum at most $m$.  We will also use Stanley's theory of $\fP$-partitions, where $\fP$ is a poset, to show that for certain $u$ and fixed $n$, these numbers are polynomials in $m$.

The next result shows, surprisingly, that when $|u|=1$ the value of $c_{n,m}(u)$ depends only on the relative sizes of $u$ and $m$ and not on their specific values.
To state it, we will use the {\em Kronecker delta function}
$$
\de_{n,m} = \case{1}{if $n=m$,}{0}{else.}
$$
\bth
\label{cnm1}
If $|u|=1$ then
$$
c_{n,m}(u)=\case{c_{n,m}(1)}{if $u\le m$,}{\de_{n,0}}{if $u>m$.}
$$
\eth
\bprf
Throughout this proof we consider $m$ as a fixed bound.
First consider the case $u>m$.  Clearly if $n=0$ then the empty word commutes with $u$, so 
$c_{0,m}(u)=1$.  On the other hand, if $n>0$ and $w\in C(u)$ then, by Theorem~\ref{|u|=1:col}, each column of $P=P(w)$ contains a $u$.  But $u>m$, the maximum value that can be used in $P$.  So $P$ does not exist and $c_{0,m}(u)=0$, finishing this case.

Now suppose $u\le m$.  It suffices to show that, for all $1\le u <m$, there is a bijection between $C(u)$ and $C(u+1)$.
Applying RSK, we get a bijection between the $w\in C(u)$ and all pairs $(P,T)$ where $P$ is a semistandard Young tableau satisfying Theorem~\ref{|u|=1:col} and $T$ is a standard Young tableau of the same shape.  So it suffices to give a shape-preserving bijection 
$\phi:\cP(u)\ra\cP(u+1)$ where
$$
\cP(u) = 
\{P \mid \text{ $P$ is an SSYT with entries at most $m$ and a $u$ in every column}\}.
$$

In fact, for $\phi$ we can just use the restriction to $\cP(u)$ of the Bender-Knuth involution interchanging $u$ and $u+1$~\cite{BK:epp}.  So,  every column of $P\in\cP(u)$ with both a $u$ and a $u+1$ remains fixed.
Any other $u$'s and $(u+1)$'s are considered free and $\phi$ interchanges the multiplicities of the free $u$'s and $(u+1)$'s in each row.

We must show $\phi$ is well defined in that $P'=\phi(P)\in\cP(u+1)$.  That is, we must show that every column of $P'$ contains a $u+1$.  This is clear for columns containing both $u$ and $u+1$.  Every other column must contain a free $u$ by definition of $\cP(u)$, and thus there are no columns containing a free $u+1$.
So these columns will all have a $u$ replaced by a $u+1$ which implies $P'\in\cP(u+1)$.  Similarly, one show that if 
$P'\in\cP(u+1)$ then applying $\phi$, thought of as a function on the set $\cP(u)\uplus\cP(u+1)$, gives
$P=\phi(P')\in\cP(u)$.  So $\phi$ is an involution and is bijective from $\cP(u)$ to $\cP(u+1)$.
\eprf

Combining the previous result with Theorem~\ref{C(a^k)}, we immediately get the following.
\bco
If  $a,k\in\bbP$ then

\vs{10pt}

\eqqed{
c_{n,m}(a^k)=\case{c_{n,m}(1)}{if $a\le m$,}{\de_{n,0}}{if $a>m$.}
}
\eco

We can now give an exact value of $c_{n,2}(a)$ for $a=1$ or $2$.
\bco
For $n\ge0$ we have
$$
c_{n,2}(1)= c_{n,2}(2) = \binom{n}{\fl{n/2}}.
$$
\eco
\bprf
By Theorem~\ref{cnm1} it suffices to prove that $c_{n,2}(1)$ is given by the corresponding central binomial coefficient.  And by Corollary~\ref{C(1)} (d), it suffices to show that if $Y_n$ is the set of Yamanouchi words in  in $[2]^n$ then
$$
\# Y_n =  \binom{n}{\fl{n/2}}.
$$
We will do this by induction on $n$ where the base case is easy to check.  So assume that the previous equation holds, and take $w=w_1 w_2\ldots w_n\in Y_n$.  Consider the set $Z_{n+1}$ obtained by prepending a $1$ or a $2$ to every element of $Y_n$.  Clearly those elements of $Z_{n+1}$ starting with a $1$ are still Yamanouchi.  There are now two cases depending on the parity of $n$.

If $n=2k+1$ is odd then any $w\in Y_n$ has at least $k+1$ ones and at most $k$ twos.  It follows that the elements of $Z_{n+1}$ ending in a $2$ are also Yamanouchi.  So, using induction,
$$
\# Y_{n+1} = \# Z_{n+1} = 2\cdot \#Y_n = 2 \binom{n}{\fl{n/2}} = 2 \binom{2k+1}{k}=\binom{2k+2}{k+1}
=\binom{n+1}{\fl{(n+1)/2}}.
$$

If $n=2k$ then the elements of $Z_{n+1}$ which are not Yamanouchi are exactly those obtained by prepending a $2$ to an element of $Y_n$ with $k$ ones and $k$ twos.  But it is well known that such sequences are counted by the $k$th Catalan number.  So, by induction again,
$$
\# Y_{n+1} =  2 \binom{n}{\fl{n/2}}-\frac{1}{k+1}\binom{2k}{k}
=\left(2 -\frac{1}{k+1}\right) \binom{2k}{k}
=\binom{2k+1}{k}
=\binom{n+1}{\fl{(n+1)/2}}
$$
finishing the proof.
\eprf

We  now show that for various $u$ and fixed $n$, 
the quantity $c_{n,m}(u)$ is a polynomial in $m$ and investigate its properties.  
We will use Stanley's theory of $\fP$-partitions where $\fP$ is a poset~\cite{sta:osp}.  For more information about this method see~\cite[Section 7.4]{sag:aoc} or~\cite[Section 3.15]{sta:ec1}.  Let $(\fP,\lte)$ be a poset on $[n]$.
Note the use of $\lte$ to differentiate the partial order in $\fP$ from the total order $\le$ on integers.
A {\em $\fP$-partition} is a function $f:\fP\ra\bbN$ satisfying
\ben
\item $i\lte j$ implies $f(i)\ge f(j)$, and
\item $i\lte j$ and $i>j$ implies $f(i)>f(j)$
\een
We let
$$
\Par_m \fP = \{f:\fP\ra[m] \mid 
\text{ $f$ is a $\fP$-partition}\}.
$$

Now suppose $\la=(\la_1,\la_2,\ldots,\la_k)$ is a partition of $n$, written $\la\ptn n$.  Partially order the cells of $\la$ reverse component-wise so that 
$(i,j)\lte(i',j')$ whenever $i\ge i'$ and $j\ge j'$.  Finally, number the cells of $\la$ with $[n]$
by numbering the first row of the Young diagram from right-to-left with $1,2,\ldots,\la_1$, then the next row right-to-left with
$\la_1+1,\la_1+2,\ldots,\la_1+\la_2$, and so forth.  Transferring this labeling to the poset constructed from $\la$ we obtain a poset $\fP_\la$.  It should be clear from the definitions that there is a bijection between the semistandard Young tableaux $P$ of shape $\la$ with maximum at most $m$ and $\Par_{m-1} \fP_\la$ obtained by subtracting one from every element of $P$.

We now describe the generating function $\sum_{m\ge0} |\Par_m\fP|\ x^m$.  If $\fP$ is a poset on $[n]$ then a  {\em linear extension} of $\fP$ is a permutation  $\pi$ in the symmetric group $\fS_n$ such that $i\lte j$  in 
$\fP$ implies $i$ is to the left of $j$ in $\pi$.  We let
$$
\cL(\fP) = \{\pi \mid \text{ $\pi$ is a linear extension of $\fP$}\}.
$$
Any $\pi=\pi_1\pi_2\ldots\pi_n\in\fS_n$ has {\em descent number}
$$
\des\pi = \#\{ i \mid \pi_i>\pi_{i+1}\}.
$$
We now have all the ingredients to state the $\fP$-partition result we will need.
\bth[\cite{sta:osp}]
\label{Par_m}
For any poset $\fP$ on $[n]$ we have

\vs{10pt}

\eqqed{
\sum_{m\ge0} |\Par_m\fP|\ x^m = \frac{\sum_{\pi\in\cL(\fP)} x^{\des\pi}}{(1-x)^{n+1}}.
}
\eth

We will now give a general lemma which will permit us to apply the previous theorem to centralizer sets.
\bth
\label{poly}
Let $u$ be a word and $r$ be a positive integer.  Suppose that  if $P=P(w)$ for $w\in C(u)$ then 
    \ben
    \item[(a)] the first $r$ rows of $P$ only contain elements which are at most $r$,
    \item[(b)] the remaining rows of $P$ can be any SSYT with elements greater than $r$.
    \een
Then for fixed $n\ge r$ and all $m\ge n$ we have that $c_{n,m}(u)$ is a polynomial in $m$ of degree $n-r$ with leading coefficient $1/(n-r)!$.
\eth
\bprf
We first note that $C(u)$ is a union of Knuth equivalence classes.  For suppose that $w\in C(u)$ and $w'\equiv w$.  Then
$$
u w' \equiv uw \equiv wu \equiv w'u
$$
so that $w'\in C(u)$.

From what we have just proved, applying RSK gives a bijection between the $w\in C(u)\cap[m]^n$ and pairs $(P,T)$ where $P$ runs over all possible seminstandard $P(w)$ and $T$ runs over all  standard Young tableaux of the same shape.  
 Suppose the mutual shape is  $\la=(\la_1,\la_2,\ldots,\la_l)\ptn n$.
 Let $g_m^\la$ denote the number of such semistandard tableaux $P$, and $f^\la$ the number of standard ones.  We have shown that
\beq
\label{gf}
c_{n,m}(u) = \sum_{\la\ptn n} g_m^\la\ f^\la.
\eeq
Since $f^\la$ is a constant for any $\la$, it suffices to show that $g_m^\la$ is a polynomial of degree at most $n-r$, and that there is only one $\la$ where it has degree $n-r$ in which case the coefficient of $m^{n-r}$ is
$1/(n-r)!$.

Given $\la$ as in the previous paragraph, we let
$\la'=(\la_1,\ldots,\la_r)$ and
$\la''=(\la_{r+1}, \ldots,\la_l)$ where $\la''\ptn n''$ for some $n''$.  
By property (a), the number of possible restrictions of $P$ to the first $r$ rows is constant with respect to $m$.  
Let $P''$ be the restriction of $P$ to $\la''$ and note that this semistandard tableau must have entries from $\{r+1,r+2,\ldots,m\}$.
Using  part (b)  and the discussion after the definion of a $\fP$-partition, there is a bijection between these $P''$ and $\fP_{\la''}$-partitions with parts at most $m-r-1$.
  It follows from Theorem~\ref{Par_m} that $g_m^\la$ is, up to a constant, the coefficient of $x^{m-r-1}$ in $(\sum_\pi x^{\des\pi})/(1-x)^{n''+1}$ where the sum is over all linear extensions of $\fP_{\la''}$.  Now 
$$
\frac{1}{(1-x)^{n''+1}}
= \sum_{m\ge0} \binom{m+n''}{n''} x^m
$$
and $\binom{m+n''}{n''}$ is a polynomial in $m$ of degree $n''$ with positive leading coefficient.
Since the numerator of the fraction in Theorem~\ref{Par_m} has positive coefficients, there will be no cancellation of leading terms so that $g_m^\la$ is also a polynomial in $m$ of degree $n''$.

To finish, note that there will only be one term of~\eqref{gf} which is of maximal degree, and that will be when $\la=(1^n)$ so that $n''=n-r$.  Furthermore, counting the tableaux involved, that term will be
$$
g_m^{(1^n)} f^{(1^n)} = \binom{m-r}{n-r}
$$
which has leading coefficient $1/(n-r)!$ as desired.
\eprf

We can now show that for fixed $n$ we have $c_{n,m}(u)$ is a polynomial in $m$ for various $u$.
\bco
The following are polynomials in $m$ for $n$ fixed and $m\ge n$.
\ben
\item[(a)] If  $n\ge 1$ then  $c_{n,m}(1)$ is a polynomial in $m$ of degree $n-1$ with leading coefficient $1/(n-1)!$.
\item[(b)] If  $n\ge 2$ then  $c_{n,m}(12)$ is a polynomial in $m$ of degree $n-2$ with leading coefficient $1/(n-2)!$.
\item[(c)] If $n\ge k$ then $c_{n,m}(k(k-1)\ldots1)$ is a polynomial in $m$ of degree $n-k$ with leading coefficient $1/(n-k)!$.
\een
\eco
\bprf
(a)  It suffices to show that the hypotheses of the previous theorem are are satisfied with $r=1$.  We get condition (a) immediately from Theorem~\ref{|u|=1:row} (a) when $u=1$.  And we see that when $u=1$ in this Theorem, condition (b) is vacuous since both sides are zero because of the row and column restrictions on a semistandard Young tableau.  Thus Theorem~\ref{poly} (b) is also true and we are done with this case.

\medskip

(b)  This is similar to (a), just using Theorem~\ref{C(12):col} in place of Theorem~\ref{|u|=1:row}.  Note that using both conditions in the former implies that the entries of $P(w)$ in the first two rows can only be $1$'s and $2$'s.  And there is no restriction on other rows.  So Theorem~\ref{poly} is satisfied with $r=2$.

\medskip

(c)  Given the similarity to (a) and (b), the reader should now be able to fill in the details using Theorem~\ref{C(m...1)}.
\eprf

One can also use the method of proof in Theorem~\ref{poly} to actually compute the polynomials involved.  We illustrate this with our next result.

\bth
\label{c_nm(1)}
Suppose $n$ is fixed and $m\ge n$.  Then we have the following polynomial expansions.
\begin{align*}
    c_{1,m}(1)&=1,\\
    c_{2,m}(1)&={m\choose 1},\\
    c_{3,m}(1)&={m\choose 1}+{m\choose 2},\\
    c_{4,m}(1)&={m\choose 1}+4{m\choose 2}+{m\choose 3},\\
    c_{5,m}(1)&={m\choose 1}+8{m\choose 2}+13{m\choose 3}+{m\choose 4},\\
    c_{6,m}(1)&={m\choose 1}+18{m\choose 2}+48{m\choose 3}+41{m\choose 4}+{m\choose 5},\\
    c_{7,m}(1)&={m\choose 1}+33{m\choose 2}+178{m\choose 3}+262{m\choose 4}+131{m\choose 5}+{m\choose 6},\\
    c_{8,m}(1)&={m\choose 1}+68{m\choose2}+549{m\choose3}+1480{m\choose4}+1405{m\choose5}+428{m\choose6}+{m\choose7}.
\end{align*}
\eth
\bprf
 We will illustrate the method by computing the term in~\eqref{gf} for $u=1$ when $\la=(2,2,1)$, leaving the details for the other partitions of integers up to $8$  to the reader.  We will continue to use the notation in the proof of Theorem~\ref{poly}.   In the case under consideration we have $\la''=(2,1)$ so that
$$
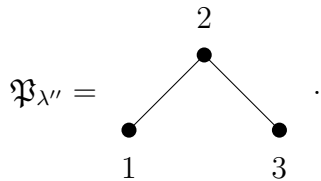
\begin{tikzpicture}
\draw(-2,.5) node{$\fP_{\la''} =$};
\fill(0,1) circle(.1);
\fill(-1,0) circle(.1);
\fill(1,0) circle(.1);
\draw (-1,0)--(0,1)--(1,0);
\draw(0,1.5) node{$2$};
\draw(-1,-.5) node{$1$};
\draw(1,-.5) node{$3$};
\draw(1.5,.5) node{.};
\end{tikzpicture}
$$
Now
$$
\cL(\fP_{\la''}) =\{ 132, 312\}
$$
so that
$$
\sum_{m\ge0} |\Par_m\fP_{\la''}|\ x^m = \frac{2x}{(1-x)^4}.
$$
The tableaux counted by $g_{5,m}^{(2,2,1)}$ correspond to $\fP_{\la''}$ partitions using entries $2,3,\ldots,m$. 
Normalizing, we want entries $0,1,\ldots,m-2$ so that we need to take the coefficient of $x^{m-2}$ in the above generating function.  This gives
$$
g_{5,m}^{(2,2,1)} =  2\binom{m}{3}.
$$
Multiplying by $f^{(2,2,1)}= 5$ then gives the desired term in~\eqref{gf}.  Finally, adding in the contributions from the other $\la\ptn 5$ results in the given expansion for $c_{5,n}(1)$.
\eprf

%%%%%%%%%%%%%%%%%%%%%%%%%%%%%%%%
\section{Open problems and conjectures}
\label{opc}
%%%%%%%%%%%%%%%%%%%%%%%%%%%%%%%%

Although we have begun the study of the centralizer $C(u)$, we believe that there are many more interesting results to be found.  Here we collect a few avenues for future research.

The reader will have noticed how useful Lemma~\ref{max<=m} was in proving various of the characterizations of $C(u)$ in Sections~\ref{csu} and~\ref{clu}.  We believe that an even stronger result is true.

\begin{conj}
       \label{maxRi}
    Given $u$, let $m=\max u$ and $\ell$ be the number of rows of $P(u)$. Suppose that $w\in C(u)$ and  that $P(w)$ has rows $R_i$ for $i\ge1$. Then
    $$
        \max R_i\leq m
    $$ for $1\leq i\leq \ell$. 
\end{conj}

To see why this conjecture implies Lemma~\ref{max<=m}, note that the existence of a subsequence of the form
$m,m-1,\ldots,m-k+1$ in $u$ implies that $\ell\ge k$ by an argument like that in the proof of the lemma.  So if $\max R_i \le m$ for $i\in[\ell]$ then certainly the inequality is true for $i\in[k]$. We have verified Conjecture~\ref{maxRi} computationally\footnote{The code used to verify the conjectures in this section can be found at \url{https://github.com/wilsoa/Centralizers-in-the-Plactic-Monoid}.} using Sage Math \cite{sagemath} for $u\in[m]^n$ and $w\in[6]^l$ where $m+n\leq 10$ and $2\leq l\leq 6$.

In Lemma~\ref{u^k} we noted that for any $u$ and $k\ge1$ we always have $C(u)\sbe C(u^k)$.
But in the particular case when $|u|=1$ we have $C(u)=C(u^k)$ for all $k\ge1$ by Theorem~\ref{C(a^k)}.  We conjecture that a similar stability phenomenon holds more generally.
\bcon\label{stability_conj}
Suppose $u\in\bbP^*$.
\ben
\item[(a)]  There is a $K\in\bbP$ such that for $k\ge K$ we have
$$
C(u^k)\sbe C(u^{k+1}).
$$
\item[(b)]  There is a $L\in\bbP$ such that for $k\ge L$ we have
$$
C(u^k)= C(u^{k+1}).
$$
\een
\econ

We have verified Conjecture~\ref{stability_conj}(a) computationally using Sage Math \cite{sagemath} for $u\in[m]^n$ and $w\in[5]^l$ where $m+n\leq 10$ and $2\leq l\leq 6$. Note that except in the particular case that $u=12345$ where $K=3$, for all other words $u$ checked, we can take $K=1$. In support of Conjecture~\ref{stability_conj}(b), the containments verified under these conditions become equalities for $k\geq4$.

The expansions of $c_{n,m}(1)$ for $n\le 8$ in Theorem~\ref{c_nm(1)} have some remarkable properties which we conjecture hold in general.  Call a sequence $a_0,a_1,\ldots,a_n$ of real numbers {\em unimodal} if there is some index $k$ such that
$$
a_0 \le a_1 \le \ldots \le a_k \ge a_{k+1} \ge \ldots a_n.
$$
The sequence is said to be {\em log-concave} if, for all $0<i<n$, 
$$
a_i^2\ge a_{i-1} a_{i+1}.
$$
Unimodal and log-concave sequences abound in combinatorics, algebra, and geometry. 
See the survey articles of Stanley~\cite{sta:lcu}, Brenti~\cite{bre:lcu}, or Br\"and\'en~\cite{bra:ulc} for more information.
It is well known that, for positive sequences, log-concavity implies unimodality.
\bcon
Fix $n$ and write
$$
c_{n,m}(1) = \sum_{k=0}^{n-1} a_k \binom{m}{k}
$$
for certain scalars $a_k$  (depending on $n$).  We have the following
\ben
\item[(a)] $a_0=0$, $a_1=1$.
\item[(b)] $a_k\in\bbP$ for all $k\in[n-1]$.
\item [(c)] The sequence $a_1, a_2,\ldots, a_{n-1}$ is log-concave and hence (assuming (b)) unimodal.
\item[(d)] The sequence $a_1, a_2,\ldots, a_{n-1}$  has maximum at $k=\ce{n/2}$.  
\een
\econ

It is well known that applying  symmetries of the square to a permutation $\pi$ (viewed as a permutation matrix) does interesting things to the output tableaux under RSK.  One of these also seems to play nicely with the centralizer set.  If 
 $w=w_1 w_2\ldots w_k$ has $\max w\le m$ then define its {\em $m$-reverse complement} to be the word
 $$
\RC_m(w) = (m-w_k+1) (m- w_{k-1}+1)\ldots(m-w_1+1).
 $$
Note the dependence on the choice of $m$, not just on $w$.  For example
$$
\RC_4(31122) = 33442.
$$
To extend this operation to an SSYT, $T$, let 
$\rw(T)$ be the row word of $T$ obtained by reading the rows of $T$ from left to right starting with the bottom row and moving up.  It is well known that 
$$
P(\rw T) = T.
$$
Now if $\max T\le m$ we define the {\em $m$-evacuation} of $T$ to be the composition
$$
\ep_m(T)= P\circ \RC_m \circ \rw(T).
$$
When $T$ is a standard Young tableau with maximum entry $m$, the map $\epsilon_m$ corresponds to 
Sch\"utzerberger's evacuation operation, see~\cite{sch:qrc} or~\cite[A1.2.10]{sta:ec2}.
\ble
\label{shape}
If $T$ is a SSYT with $\max T\le m$, then $T$ and $\ep_m(T)$ have the same shape.
\ele
\bprf
Let $T$ have shape $\la=(\la_1,\la_2,\ldots,\la_l)$.  Then by Greene's extension of Schensted's Theorem~\cite{gre:est,sch:lid} we have that, for all $i\ge 1$, the sum $\la_1+\la_2+\cdots+\la_i$ is the cardinality of the largest subword of $w=\rw T$ which can be written as a disjoint union of $i$ weakly increasing subwords.
So, to prove the lemma, it suffices to give a bijection between weakly increasing subwords of $w$ and those of $\RC_m(w)$ which preserves disjointness.
But $v$ is a weakly increasing subword of $w$ if and only if $\RC_m(v)$ is a weakly increasing subword of $\RC_m(w)$, so we are done.
\eprf

We now wish to describe a conjectural bijection between the elements of $C(u)$ and
those of $C(\RC_m(u))$ for any $m\ge \max u$.
Given an array $A$ we let $A_{\leq m}$  be the subarray consisting of the elements of $A$ which are at most $m$ and similarly for $A_{>m}$.
If $T$ is an SSYT then let $\tau_m(T)$ be the result of replacing $T_{\leq m}$ with $\epsilon_m(T_{\leq m})$ and leaving $T_{>m}$ unchanged.  Note that the previous lemma makes this replacement well defined since the two tableaux involved have the same shape.  We similarly extend the map $\RC_m$ to all words $w$ by letting $\RC_m(w)$ be the word obtained by replacing $w_{\le m}$ with its $m$-reverse complement and leaving the elements of $w_{>m}$ unchanged.  

As noted in the proof of Theorem~\ref{poly}: for any $u$, the centralizer $C(u)$ is a union of Knuth equivalence classes.  It follows that to describe $C(u)$ it suffices to describe the set $P(C(u))$ of all $P(w)$ for $w\in C(u)$.
Our conjecture is as follows.
\bcon
    \label{RC_conj}
    If $u$ is a word with $\max u\le m$ then
    $$
    P( C( \RC_m(u) ) ) = \tau_m ( P( C(u) ) ).    
    $$
\econ

Note that, since both $\RC_m$ and $\tau_m$ are involutions, it suffices to prove only one of the two set containments implied by the conjectured equality.
We have verified this computationally 
%that $w\in C(u)$ implies that 
%$$
%\rw(\tau_m(P(w)))\in C(\RC_m(u))
%$$ 
%where $\max{u}\leq m$ 
for $u\in[m]^n$ and $w\in[6]^l$ where $m+n\leq 11$ and $2\leq l\leq 5$. 
%To see how this statement implies Conjecture~\ref{RC_conj} in these cases, we first note that $\RSK^{-1}(\tau_m(P(w)),Q(w))\equiv \rw(\tau_m(\RC_m(w)))$. Hence, the image of $C(u)$ under the map in Conjecture~\ref{RC_conj} does indeed land in $C(\RC_m(u))$. Because this map and $\RC_m$ are involutions, it is furthermore clear that the image of $C(u)$ under this map is precisely $C(\RC_m(u))$ and that the map is a bijection.

%\bibliographystyle{plain}

%\nocite{*}
%\bibliographystyle{abbrvnat}

\nocite{*}
\bibliographystyle{alpha}

\end{document}